\documentclass[10pt]{amsart}

\usepackage{amssymb, amscd, amsmath, amsthm,  graphicx, latexsym, psfrag}

\def\zz{{\bf Z}}

\newtheorem{theorem}{Theorem}
\newtheorem{proposition}[theorem]{Proposition}
\newtheorem{lemma}[theorem]{Lemma}

\newtheorem{definition}[theorem]{Definition}

\begin{document}

\title[An obstruction to slicing iterated Bing doubles]
{An obstruction to slicing iterated Bing doubles}

\author{Cornelia A. Van Cott}
\address{
University of San Francisco, San Francisco, California 94117} 
\email{cvancott@usfca.edu}
\maketitle

\begin{abstract} 
Let $K$ be a knot in $S^3$.  We study the iterated Bing doubles of $K$, giving a new proof for the following statement:  If $BD_n(K)$ is slice for some $n$, then $K$ is algebraically slice.  This result was first proved by Cha and Kim using covering link calculus.  We also use this tool, but our proof is substantially simpler and illuminates several generalizations.

\end{abstract}

\section{Introduction}
\begin{figure}
\begin{center}
\includegraphics[width=6cm]{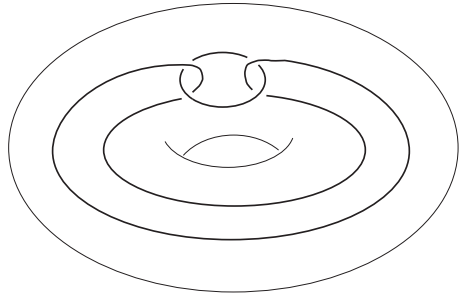}
\end{center}
\caption{A solid torus $T$ containing a link of two components. }\label{solid}
\end{figure}

The Bing double of a knot $K$ is a two component link.  The link is formed by finding a homeomorphism of the solid torus $T$ in Figure~\ref{solid} onto a tubular neighborhood of $K$, where the homeomorphism takes meridian to meridian and longitude to longitude.  The image of the two components inside $T$ under this homeomorphism is the Bing double of $K$, denoted $BD(K)$.  

\begin{figure}
\begin{center}
\includegraphics[width=13cm]{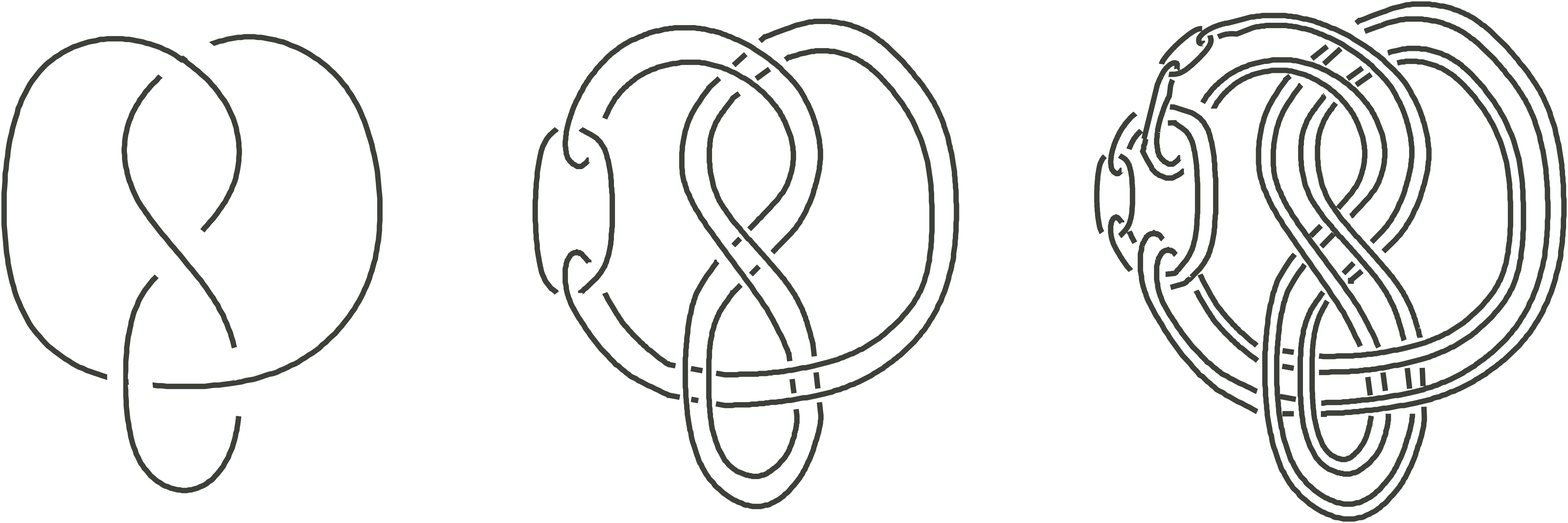}
\end{center}
\caption{The figure eight knot $K$, its Bing double $BD(K)$, and its second iterated Bing double $BD_2(K)$. }\label{example}
\end{figure}

More generally, the Bing double of a link $L$ is constructed by taking a set of homeomorphisms which map the solid torus $T$ onto a tubular neighborhood of each component of $L$, respectively, and considering union of the images of the two components inside of $T$ under each homeomorphism.  If $L$ has $n$ components, then the Bing double of $L$ has $2n$ components.  In this way, the process of Bing doubling can be iterated.  Given any knot $K$, we have a sequence of links $BD_n(K)$, called iterated Bing doubles of $K$.  See Figure~\ref{example} for an example.

Bing doubling originally arose from a 3-manifold construction of R.~H.~Bing~\cite{bing}.  Many years later, Bing doubles surfaced in the work of Freedman in an effort to solve a long-standing conjecture in 4-manifold theory~\cite{freedman}.  Freedman's work concentrated on showing that particular iterated Bing doubles were slice.  This work motivated more general study of when Bing doubles are slice.  (Recall that a link is slice in $S^3$ if the components bound disjoint topologically locally flat disks in $B^4$.)

Investigating when Bing doubles are slice is complicated by the fact that many of the traditional link invariants which obstruct sliceness vanish on all Bing doubles.  For example, the multivariable Alexander polynomial, multivariable Tristam-Levine signatures, the Arf invariant, and Milnor's invariants {\em all} vanish on $BD(K)$ for all $K$ (for an excellent summary of these results, see~\cite{cimasoni}).  Hence, subtler methods are required to obstruct sliceness.

Recent attention has focused on discovering how the slicing of the Bing double $BD(K)$ is related to slicing the associated knot $K$.  For instance, one can easily show that if a knot $K$ is slice, then the Bing double $BD(K)$ is also slice.  A natural question to ask is whether all slice Bing doubles arise in this way.  That is, if a Bing double $BD(K)$ is slice, is the knot $K$ slice?  

Progress toward proving such a statement has been incremental.  If the link $BD(K)$ is slice, the strongest conclusion which can be made at this time is that the knot $K$ must be {\em algebraically slice}, a weaker notion than sliceness which requires a particular bilinear pairing associated to a knot to be metabolic.  

The history behind this result begins with Harvey~\cite{harvey} and Teichner (unpublished) who proved that if $BD(K)$ is slice, then the integral of the Tristam-Levine signature function is zero.  Cimasoni~\cite{cimasoni} subsequently proved that if $BD(K)$ is boundary slice, then $K$ is algebraically slice.   Next, using new machinery of 3-manifold invariants, Cha~\cite{cha} strengthened the signature result of Harvey and Teichner by proving that if $BD(K)$ is slice, then the entire signature function vanishes.  Cha-Livingston-Ruberman~\cite{clr} later proved the following:\\

\begin{theorem}~\label{livingston}~\cite{clr}
If $BD(K)$ is slice, then the knot $K$ is algebraically slice.  
\end{theorem}

This result was subsequently generalized by Cha-Kim to the family of iterated Bing doubles $BD_n(K)$ as follows:

\begin{theorem}\label{taehee}~\cite{chataehee}
If $BD_n(K)$ is slice for some $n$, then $K$ is algebraically slice.
\end{theorem}

The purpose of this note is to provide a new proof for Theorem~\ref{taehee}.  The tools used will be the same as the original proof, but the argument is substantially simpler and illuminates several generalizations.

%In~\cite{}, Cha and Kim develop a technique called {\em covering link calculus} to study links.  We will use this technique here, as well.  Therefore, we begin in Section~\ref{background} by reviewing the details of the technique which are necessary to our argument.  For a more complete description, the reader is referred to Cha and Kim's work~\cite{}.  In Section~\ref{}, 

%Before outlining the proof strategy, we note that the converse of Theorem~\ref{taehee} is false.  By looking at higher order signatures, Cochran-Harvey-Leidy~\cite{chl} found that for every $n$, there exist examples of algebraically slice knots for which $BD_n(K)$ is not slice.  In addition, Cha~\cite{cha} found an infinite family of amphicheiral knots with non-slice iterated Bing doubles.  These examples of Cha were the first examples of a knot $K$ with finite concordance order such that the Bing double is not slice.  % CHECK THIS!?!?!?!??!?!?
%\section{Covering links}~\label{background}

\section{Background}
To begin, let $\zz_{(p)}$ denote $\zz$ localized at $p$.  That is, 
$$\zz_{(p)} = \{ a / b \mid a, b \in \zz, p \nmid b\}.$$
This ring $\zz_{(p)}$ surfaces when studying branched covers of $S^3$ as follows.  
\begin{lemma}~\label{primefoldcover}
Let $p$ be a prime.  The $p^k$-fold cyclic branched cover of $S^3$ branched over a knot $K$ is a $\zz_{(p)}$-homology sphere.
\end{lemma}

With this in mind, we introduce a different flavor of sliceness:\\

\begin{definition}~\label{z2slice}
A link $L$ in a $\zz_{(p)}$-homology sphere $X$  is said to be {\em $\zz_{(p)}$-slice} if there is a $\zz_{(p)}$-homology ball $Y$ with boundary $X$ such that the components of $L$ bound disjoint topologically locally flat disks in $Y$.
\end{definition}

Notice that if a link $L\subset S^3$ is slice, then it is $\zz_{(p)}$-slice for all $p$.  One of the useful results about $\zz_{(p)}$-slice links in $S^3$ is that the lifts of these links in branched covers remain $\zz_{(p)}$-slice. 

\begin{theorem}~\label{zp}
Let $p$ be any prime and $L$ a link in $S^3$.  Let $M$ be a branched cover of $S^3$ branched over one of the components of $L$.  If $L$ is $\zz_{(p)}$-slice, then the lift of $L$ in $M$ is also $\zz_{(p)}$-slice.
\end{theorem}

For a proof of the above theorem, see Cha and Kim's paper~\cite{chataehee}.  Restricting to the case when $p = 2$, we have the following result (see, for example,~\cite{clr} for a proof):

\begin{theorem}  
If $K \subset S^3$ bounds a slice disk in a $\zz_{(2)}$-homology ball $Y$, then $K$ is algebraically slice.
\end{theorem}

An analogous result does not hold if $\zz_{(2)}$ is replaced by $\zz_{(p)}$~\cite{clr}.  

\section{Proof of Theorem~\ref{taehee}}
The proof we present here is inductive in flavor.  We prove the theorem in two steps:

\begin{proposition}~\label{downone}
Let $n\geq 2$.  If $BD_{n}(K)$ is $\zz_{(2)}$-slice, then $BD_{n-1}(K)$ is $\zz_{(2)}$-slice.
\end{proposition}

\begin{proposition}\label{caseone}
If $BD_1(K)$ is $\zz_{(2)}$-slice, then $K$ is algebraically slice.
\end{proposition}

Theorem~\ref{taehee} follows from these two propositions, since $BD_n(K)$ being slice implies, in particular, that $BD_n(K)$ is $\zz_{(2)}$-slice.  

We begin with the easier step: Proposition~\ref{caseone}.  This statement is similar to the statement of Theorem~\ref{livingston}.  In fact, carefully reading the proof of Theorem~\ref{livingston} in \cite{clr}, it is clear that the argument does not use the full assumption that the link is slice.  Instead, the argument only uses that the link is $\zz_{(2)}$-slice. Hence the proposition follows.

Now to prove Proposition~\ref{downone}, we begin with the $n^{th}$ iterated Bing double of $K$, $BD_{n}(K)$, where $n \geq 2$.  By assumption, the link is $\zz_{(2)}$-slice.  Recall that this link is constructed from the $(n-1)^{th}$ iterated Bing double by replacing each individual component in $BD_{n-1}(K)$ with two components.  There are $2^{n-1}$ components in $BD_{n-1}(K)$.  Hence we can think of the link $BD_{n}(K)$ as lying within a collection of $2^{n-1}$ solid tori, each of which contains 2 components of the link.  Notice that since $n\geq 2$, these solid tori are themselves unknotted.  

\begin{figure}

\begin{center}
\psfrag{J1}{\!\!\!\footnotesize $J_1$}
\psfrag{J2}{\!\!\!\footnotesize $J_2$}
\psfrag{T}{\!\!\!\footnotesize $T$}
\includegraphics[width=8cm]{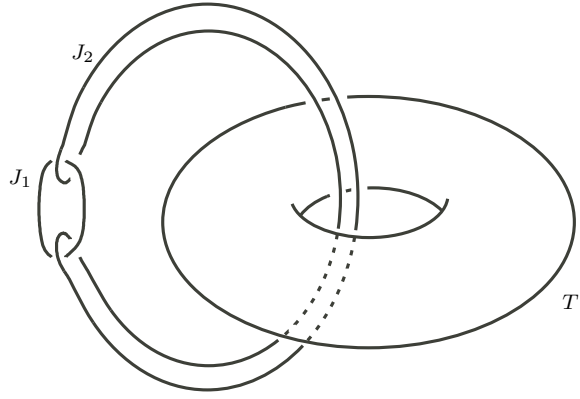}
\end{center}

\caption{The link $BD_{n}(K)$ can be drawn with two components as above, and the remaining components inside the solid torus $T$. }\label{formal}
\end{figure}

\begin{figure}
\begin{center}
\psfrag{a}{\scriptsize ${J_1}$}
\psfrag{b}{$^{J_2'}$}
\psfrag{c}{\scriptsize $T'$}
\psfrag{d}{$^{ J_2''}$}
\psfrag{e}{\scriptsize $T''$}
\includegraphics[width=13cm]{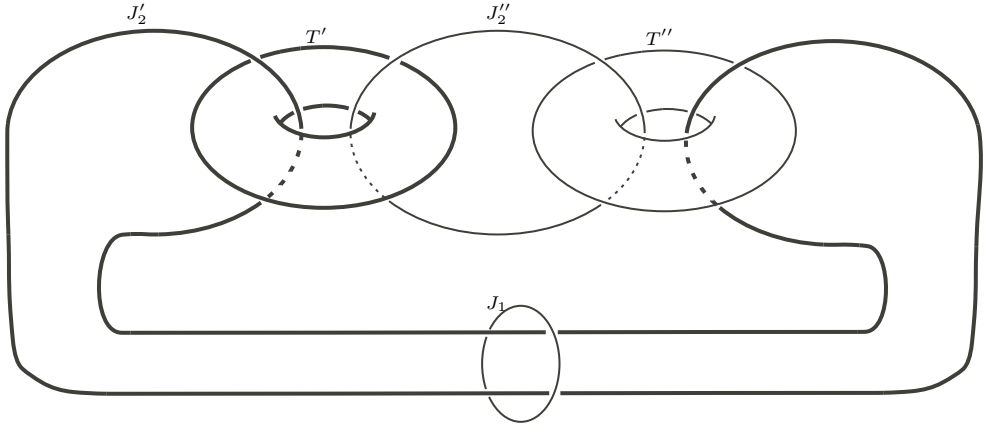}
\end{center}
\caption{The lift of $BD_{n}(K)$ pictured in Figure~\ref{formal} to the cover of $S^3$ branched over $J_1$. }\label{cover}
\end{figure}

\begin{figure}
\begin{center}
\psfrag{J}{\footnotesize $J_2'$}
\psfrag{T}{\footnotesize $T'$}
\includegraphics[width=8cm]{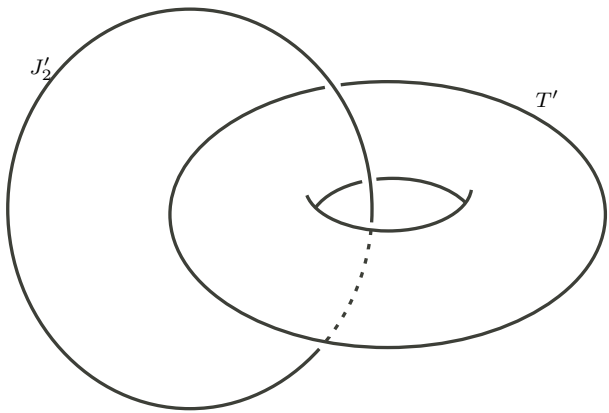}
\end{center}
\caption{A sublink $\widetilde L$ of the link in Figure~\ref{cover}.  }\label{sublink}
\end{figure}

Consider any one of these unknotted solid tori.  Denote this solid torus by $S$, and denote the two components inside $S$ by $J_1$ and $J_2$.  The link $BD_{n}(K)$ can be drawn as in Figure~\ref{formal}, where the torus $S$ is in the exterior of the torus marked $T$, and $T$ contains all the remaining tori and components of the link.  Take the $2$-fold  branched cover over the component denoted $J_1$ in Figure~\ref{formal}.  Notice that $J_1$ is unknotted, and so the branched cover of $S^3$ branched over $J_1$ is again $S^3$.  The lift of the link $BD_{n}(K)$ upstairs in the cover is denoted in Figure~\ref{cover}, using the same notation as in Figure~\ref{formal}.  Since $BD_{n}(K)$ is $\zz_{(2)}$-slice, it follows by Theorem~\ref{zp} that the lift of the link upstairs is $\zz_{(2)}$-slice, as well.  Consider the sublink $\widetilde L$ of the covering link consisting of $J_2'$ together with all of the components inside the solid torus $T'$, as highlighted in Figure~\ref{cover}.  Observe that if we forget about the other components of the covering link, the components of $\widetilde L$ can be redrawn as in Figure~\ref{sublink}.  As a sublink of a $\zz_{(2)}$-slice link, $\widetilde L$ is $\zz_{(2)}$-slice.  

The link $\widetilde L$ is identical to $BD_{n}(K)$ except the process of Bing doubling has been reversed in one of the tori.  Repeat this process on each of the remaining solid tori one by one:  take the branched cover over one of the components in the torus, look at the lift of the link upstairs, and take the sublink denoted in Figure~\ref{sublink}.  Repeating this process, the result will be a tower of $2$-fold covers, and a tower of links, each of which is $\zz_{(2)}$-slice.  At the top of the tower of links will be a link in which the process of Bing doubling has been reversed in all the solid tori.  Hence this link is $BD_{n-1}(K)$.  As a sublink of a $\zz_{(2)}$-slice link, it is also $\zz_{(2)}$-slice.  This completes the proof.

\section{Other satellite links}
Our discussion thus far has been limited to the family of Bing doubles, however the methods we have discussed can be used to obstruct other links from being slice, as well.  Bing doubles belong to a larger family of links called satellite links.  A satellite link is formed in the same way as a Bing double, except that the link inside the torus pictured in Figure~\ref{solid} is replaced by an arbitrary link $L$.   We denote the resulting link by $L(K)$.  Iterating the satellite construction, we obtain a sequence of links $L_n(K)$.

\begin{figure}
\begin{center}
\includegraphics[width=8cm]{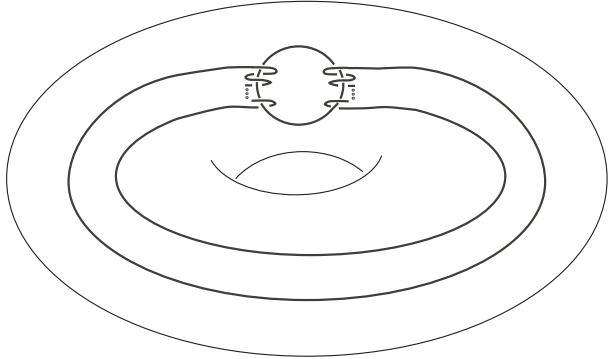}
\end{center}
\caption{A solid torus $T$ containing the link $L$, where there are $r$ full twists on either side of the small component.  Notice that if $r = 1$, then we have again the link in Figure~\ref{solid}. }\label{solid2}
\end{figure}

Our argument for finding an obstruction to slicing iterated Bing doubles generalizes to a large class of iterated satellite links and offers partial information in other cases.  In particular, let $L$ be a link of two components in a solid torus where both components of $L$ are unknotted, each component is null-homologous in the solid torus, and the linking number of the two components is zero.  (See Figures~\ref{solid} and ~\ref{solid2} for examples.)  Under these conditions, we can retool the analysis given here for Bing doubles to find an obstruction to the links $L_n(K)$ being slice.  

For example, consider the link $L$ in Figure~\ref{solid2}.  We can show that if $L_n(K)$ is slice for some $n$, then $K$ is algebraically slice.  This example, as well as other examples, will be discussed in detail in~\cite{vancott}, where we will carefully develop a general procedure from obstructing families of links of this form from being slice.
\vspace{1cm}

\bibliographystyle{plain}
\bibliography{bibliography}  

\end{document}